\documentclass{article}
\usepackage{amssymb}
\usepackage{amsmath}
\usepackage{latexsym,amssymb,a4wide}
\usepackage{latexsym,xypic}

\newtheorem{theorem}{Theorem}

\newtheorem{definition}[theorem]{Definition}

\newtheorem{lemma}[theorem]{Lemma}

\newtheorem{proposition}[theorem]{Proposition}
\newtheorem{remark}[theorem]{Remark}

\newenvironment{proof}[1][Proof]{\noindent\textbf{#1.} }{\ \rule{0.5em}{0.5em}}

\begin{document}

\date{}
\author{T.S. Kuzp\i nar\i , A. Odaba\c{s} and E.\"{O}. Uslu}
\title{On $3$-crossed modules of algebras}
\maketitle

\begin{abstract}
In this paper we define 3-crossed modules for commutative (Lie) algebras and
investigate the relation between this construction and the simplicial
algebras. Also we define the projective $3$-crossed resolution for
investigate a higher dimensional homological information \ and show the
existence of this resolution for an arbitrary $\mathbf{k}$-algebra.
\end{abstract}


\section{Introduction}

As an extension of crossed modules (Whitehead) and $2$-crossed modules
(Conduch\'{e}); Arvasi, Kuzp\i nar\i\ and Uslu in \cite{arkaus}, defined $3$%
-crossed modules as a model for homotopy 4-types. Kan, in \cite{kan}, also
proved that simplicial groups are algebraic models for homotopy types. It is
known from \cite{patron3, grand, portervascon} that simplicial algebras with
Moore complex of length $1,$ $2$ lead to crossed module and $2$-crossed
modules which are related to Koszul complex and Andre-Quillen homology
constructions for use in homotopical and homological algebra. \bigskip

PJ.L.Doncel, A.R. Grandjean and M.J.Vale in \cite{doncel} extent the $2$%
-crossed modules of groups to commutative algebras. As in indicated in \cite%
{doncel}, they defined a homology theory and obtain the relation with
Andre-Quillen homology for $n=0,1,2,$. This homology theory includes the
projective $2$-crossed \ resolution and the homotopy operator given in \cite%
{vale1}. Of course these results based on the work of T.Porter \cite%
{portervascon}, which involves the relation between Koszul complex and the
Andre-Quillen homology by means of free crossed modules of commutative
algebra. In this vein, we hope that it would be possible to generalise these
results by using commutative algebra case of higher dimensional crossed
algebraic gadgets. \bigskip

The present work involves the relation between $3$-crossed modules and
simplicial algebra without details since the most calculation are same as
group case given in \cite{arkaus}. Furthermore the work involves the
existence of projective 3-crossed resolution of a $\mathbf{k}$\textbf{-}%
algebra to obtain an higher dimensional homological information about
commutative algebras. Here the construction is a bit different from the
2-crossed resolution given in \cite{doncel} because of the number of Peiffer
liftings. At the end of the work we give the Lie algebra 3-crossed modules.

\bigskip

The main results of this work are;

\begin{enumerate}
\item Introduce the notion of $3$-crossed modules of commutative algebras
and Lie algebras;

\item Construct the passage from $3$-crossed modules of algebras to
simplicial algebras and the converse passage as an analogue result given in
\cite{arkaus};

\item Define the projective $3$-crossed resolution for investigate a higher
dimensional homological information \ and show the existence of this
resolution for an arbitrary $\mathbf{k}$-algebra which was shown for two
crossed modules in \cite{doncel}.
\end{enumerate}

\section{Preliminaries}

In this work $\mathbf{k}$ will be a fixed commutative ring with identity $1$
not equal to zero and all algebras will be commutative $\mathbf{k}$\textbf{-}%
algebras, we accept they are not required to have the identity $1$.

\subsection{\textbf{Simplicial Algebras }}

See \cite{may}, \cite{curtis} for most of the basic properties of simplicial
structures.

A simplicial algebra $\mathbf{E}$ consists of a family of algebras $\left\{
E_{n}\right\} $ together with face and degeneracy maps $d_{i}^{n}:E_{n}%
\rightarrow E_{n-1}$, $0\leq i\leq n$, $(n\neq 0)$ and $s_{i}^{n}:E_{n-1}%
\rightarrow E_{n},$ $0\leq i\leq n$, satisfying the usual simplicial
identities given in \cite{andre}, \cite{ill}. The category of simplicial
algebras will be denoted by $\mathbf{SimpAlg}$.

Let $\Delta $ denotes the category of finite ordinals. For each $k\geq 0$ we
obtain a subcategory $\Delta _{\leq k}$ determined by the objects $\left[ i%
\right] $ of $\Delta $ with $i\leq k$. A $k$-truncated simplicial algebras
is a functor from $\Delta _{\leq k}^{op}$ to $\mathbf{Alg }$ (the category
of algberas). We will denote the category of $k$-truncated simplicial
algebras by $\mathbf{Tr}_{k}\mathbf{SimpAlg }\mathfrak{.}$ By a $k$-$%
truncation$ $of$ $a$ $simplicial$ $algebra,$ we mean a $k$-truncated
simplicial algebra $\mathbf{tr}_{k}\mathbf{E}$ obtained by forgetting
dimensions of order $>k$ in a simplicial algebra $\mathbf{E.}$ Then we have
the adjoints situations

\begin{center}
$\xymatrix@R=60pt@C=60pt{\ \mathbf{SimpAlg} \ar@{->}@<2pt>[r]^-{\mathbf{tr}%
_{k}} & \mathbf{Tr}_{k} \mathbf{SimpAlg} \ar@{->}@<2pt>[l]^-{\mathbf{st}_{k}}
\ar@{->}@<2pt>[r]^-{\mathbf{cost}_{k}} & \mathbf{SimpAlg} \ar@{->}@<2pt>[l]^-%
{\mathbf{tr}_{k}} }$ 
\end{center}

\noindent where $\mathbf{st}_{k}$ and $\mathbf{cost}_{k\text{ }}$are called
the $k$-skeleton and the $k$-coskeleton functors respectivily For detailed
definitions see \cite{duskin}.

\subsection{\textbf{The Moore Complex}.}

The Moore complex $\mathbf{NE}$\ of a simplicial algebra $\mathbf{E}$\ is
defined to be the normal chain complex $(\mathbf{NE,\partial })$ with%
\begin{equation*}
NE_{n}=\bigcap\limits_{i=0}^{n-1}kerd_{i}
\end{equation*}%
and with differential $\partial _{n}:NE_{n}\rightarrow NE_{n-1}$ induced
from $d_{n}$ by restriction.

We say that the Moore complex $\mathbf{NE}$ of a simplicial algebra $\mathbf{%
E}$ is of \textit{length k} if $\mathbf{NE}_{n}=0$ for all $n\geq k+1$. We
denote the category of simplicial algebras with Moore complex of length $k$
by $\mathbf{SimpAlg}_{\leq k}.$

The Moore complex, $\mathbf{NE}$, carries a hypercrossed complex structure
(see Carrasco \cite{c} ) from which $\mathbf{E}$ can be rebuilt. Now we will
have a look to this construction slightly. The details can be found in \cite%
{c}.

\subsection{\textbf{The Poset of Surjective Maps}}

The following notation and terminology is derived from \cite{cc}.

For the ordered set $[n]=\{0<1<\dots <n\}$, let $\alpha
_{i}^{n}:[n+1]\rightarrow \lbrack n]$ be the increasing surjective map given
by;
\begin{equation*}
\alpha _{i}^{n}(j)=\left\{
\begin{array}{ll}
j & \text{if }j\leq i, \\
j-1 & \text{if }j>i.%
\end{array}%
\right.
\end{equation*}%
Let $S(n,n-r)$ be the set of all monotone increasing surjective maps from $%
[n]$ to $[n-r]$. This can be generated from the various $\alpha _{i}^{n}$ by
composition. The composition of these generating maps is subject to the
following rule: $\alpha _{j}\alpha _{i}=\alpha _{i-1}\alpha _{j},j<i$. This
implies that every element $\alpha \in S(n,n-r)$ has a unique expression as $%
\alpha =\alpha _{i_{1}}\circ \alpha _{i_{2}}\circ \dots \circ \alpha
_{i_{r}} $ with $0\leq i_{1}<i_{2}<\dots <i_{r}\leq n-1$, where the indices $%
i_{k}$ are the elements of $[n]$ such that $\{i_{1},\dots
,i_{r}\}=\{i:\alpha (i)=\alpha (i+1)\}$. We thus can identify $S(n,n-r)$
with the set $\{(i_{r},\dots ,i_{1}):0\leq i_{1}<i_{2}<\dots <i_{r}\leq
n-1\} $. In particular, the single element of $S(n,n)$, defined by the
identity map on $[n]$, corresponds to the empty 0-tuple ( ) denoted by $%
\emptyset _{n} $. Similarly the only element of $S(n,0)$ is $(n-1,n-2,\dots
,0)$. For all $n\geq 0$, let
\begin{equation*}
S(n)=\bigcup_{0\leq r\leq n}S(n,n-r).
\end{equation*}%
We say that $\alpha =(i_{r},\dots ,i_{1})<\beta =(j_{s},\dots ,j_{1})$ in $%
S(n)$ if $i_{1}=j_{1},\dots ,i_{k}=j_{k}$ but $i_{k+1}>j_{k+1},$ $(k\geq 0)$
or if $i_{1}=j_{1},\dots ,i_{r}=j_{r}$ and $r<s$. This makes $S(n)$ an
ordered set. For example

\begin{eqnarray*}
S(2) &=&\{\phi _{2}<(1)<(0)<(1,0)\} \\
S(3) &=&\{\phi _{3}<(2)<(1)<(2,1)<(0)<(2,0)<(1,0)<(2,1,0)\} \\
S(4) &=&\{\phi _{4}<(3)<(2)<(3,2)<(1)<(3,1)<(2,1)<(3,2,1)<(0)<(3,0)<(2,0) \\
&<&(3,2,0)<(1,0)<(3,1,0)<(2,1,0)<(3,2,1,0)\}
\end{eqnarray*}

\subsection{The Semidirect Decomposition of a Simplicial Algebra}

The fundamental idea behind this can be found in Conduch\'{e} \cite{conduche}%
. A detailed investigation of this for the case of simplicial groups is
given in Carrasco and Cegarra \cite{cc}.The algebra case of the structure is
also given in \cite{c}.

\begin{proposition}
If \textbf{E} is a simplicial algebra, then for any $n\geq 0$%
\begin{equation*}
\begin{array}{lll}
E_n & \cong & (\ldots (NE_n \rtimes s_{n-1}NE_{n-1})\rtimes \ldots \rtimes
s_{n-2}\ldots s_0NE_1)\rtimes \\
&  & \qquad (\ldots (s_{n-2}NE_{n-1}\rtimes s_{n-1}s_{n-2}NE_{n-2})\rtimes
\ldots \rtimes s_{n-1}s_{n-2}\dots s_0NE_0). ~~%
\end{array}%
\end{equation*}
\end{proposition}

\begin{proof}
This is by repeated use of the following lemma.
\end{proof}

\begin{lemma}
Let \textbf{E} be a simplicial algebra. Then $E_{n}$ can be decomposed as a
semidirect product:
\begin{equation*}
E_{n}\cong \mathrm{ker}d_{n}^{n}\rtimes s_{n-1}^{n-1}(E_{n-1}).
\end{equation*}
\end{lemma}

The bracketing and the order of terms in this multiple semidirect product
are generated by the sequence:

\begin{equation*}
\begin{array}{lll}
E_{1} & \cong & NE_{1}\rtimes s_{0}NE_{0} \\
E_{2} & \cong & (NE_{2}\rtimes s_{1}NE_{1})\rtimes (s_{0}NE_{1}\rtimes
s_{1}s_{0}NE_{0}) \\
E_{3} & \cong & ((NE_{3}\rtimes s_{2}NE_{2})\rtimes (s_{1}NE_{2}\rtimes
s_{2}s_{1}NE_{1}))\rtimes \\
&  & \qquad \qquad \qquad ((s_{0}NE_{2}\rtimes s_{2}s_{0}NE_{1})\rtimes
(s_{1}s_{0}NE_{1}\rtimes s_{2}s_{1}s_{0}NE_{0})).%
\end{array}%
\end{equation*}%
and%
\begin{equation*}
\begin{array}{lll}
E_{4} & \cong & (((NE_{4}\rtimes s_{3}NE_{3})\rtimes (s_{2}NE_{3}\rtimes
s_{3}s_{2}NE_{2}))\rtimes \\
&  & \qquad \ ((s_{1}NE_{3}\rtimes s_{3}s_{1}NE_{2})\rtimes
(s_{2}s_{1}NE_{2}\rtimes s_{3}s_{2}s_{1}NE_{1})))\rtimes \\
&  & \qquad \qquad s_{0}(\text{decomposition of }E_{3}).%
\end{array}%
\end{equation*}

Note that the term corresponding to $\alpha =(i_{r},\ldots ,i_{1})\in S(n)$
is%
\begin{equation*}
s_{\alpha }(NE_{n-\#\alpha })=s_{i_{r}...i_{1}}(NE_{n-\#\alpha
})=s_{i_{r}}...s_{i_{1}}(NE_{n-\#\alpha }),
\end{equation*}%
where $\#\alpha =r.$ Hence any element $x\in E_{n}$ can be written in the
form%
\begin{equation*}
x=y+\sum\limits_{\alpha \in S(n){\backslash }\left\{ \emptyset _{n}\right\}
}s_{\alpha }(x_{\alpha })\ \text{with }y\in NE_{n}\ \text{and }x_{\alpha
}\in NE_{n-\#\alpha }.
\end{equation*}

\subsection{\textbf{Hypercrossed Complex Pairings}}

In the following we recall from \cite{patron3} hypercrossed complex pairings
for commutative algebras. The fundamental idea behind this can be found in
Carrasco and Cegarra (cf. \cite{cc}). The construction depends on a variety
of sources, mainly Conduch\'{e} \cite{conduche}, Z. Arvasi and T. Porter,
\cite{patron3}. Define a set $P(n)$ consisting of pairs of elements $(\alpha
,\beta )$ from $S(n)$ with $\alpha \cap \beta =\emptyset $ and $\beta
<\alpha $ , with respect to lexicographic ordering in $S(n)$ where $\alpha
=(i_{r},\dots ,i_{1}),\beta =(j_{s},\dots ,j_{1})\in S(n)$. The pairings
that we will need,
\begin{equation*}
\{C_{\alpha ,\beta }:NE_{n-\sharp \alpha }\otimes NE_{n-\sharp \beta
}\rightarrow NE_{n}:(\alpha ,\beta )\in P(n),n\geq 0\}
\end{equation*}%
are given as composites by the diagram

\begin{center}
$\xymatrix@R=40pt@C=60pt{\ NE_{n-\#\alpha}\otimes NE_{n-\#\beta} \ar[d]%
_{s_{\alpha}\otimes s_{\beta}} \ar[r]^-{C_{\alpha ,\beta}} & NE_n \\
E_n \otimes E_n \ar[r]_{\mu} & E_n \ar[u]_{p} }$ 
\end{center}

where $s_{\alpha }=s_{i_{r}},\dots ,s_{i_{1}}:NE_{n-\sharp \alpha
}\rightarrow E_{n},$\quad $s_{\beta }=s_{j_{s}},\dots
,s_{j_{1}}:NE_{n-\sharp \beta }\rightarrow E_{n},$ \newline
$p:E_{n}\rightarrow NE_{n}$ is defined by composite projections $%
p(x)=p_{n-1}\dots p_{0}(x),$ where \newline
$p_{j}(z)=zs_{j}d_{j}(z)^{-1}$ with $j=0,1,\dots ,n-1.$ $\mu :E_{n}\otimes
E_{n}\rightarrow E_{n}$ is given by multiplication map and $\sharp \alpha $
is the number of the elements in the set of $\alpha ,$ similarly for $\sharp
\beta .$ Thus%
\begin{eqnarray*}
C_{\alpha ,\beta }(x_{\alpha }\otimes y_{\beta }) &=&p\mu (s_{\alpha
}\otimes s_{\beta })(x_{\alpha }\otimes y_{\beta }) \\
&=&p(s_{\alpha }(x_{\alpha })\otimes s_{\beta }(y_{\beta })) \\
&=&(1-s_{n-1}d_{n-1})\dots (1-s_{0}d_{0})(s_{\alpha }(x_{\alpha })s_{\beta
}(y_{\beta }))
\end{eqnarray*}

Let $I_{n}$ be the ideal in $E_{n}$ generated by elements of the form
\begin{equation*}
C_{\alpha ,\beta }(x_{\alpha }\otimes y_{\beta })
\end{equation*}
where $x_{\alpha }\in NE_{n-\sharp \alpha }$ and $y_{\beta }\in NE_{n-\sharp
\beta }.$

We illustrate this for $n=3$ and $n=4$ as follows:

For $n=3$, the possible Peiffer pairings are the following

\begin{center}
$C_{(1,0)(2)}$, $C_{(2,0)(1)}$, $C_{(0)(2,1)}$, $C_{(2)(0)}$, $C_{(2)(1)}$, $%
C_{(1)(0)}$
\end{center}

For all $x_{1}\in NE_{1},y_{2}\in NE_{2},$ the corresponding generators of $%
I_{3}$ are:
\begin{align*}
C_{(1,0)(2)}(x_{1}\otimes y_{2})&
=(s_{1}s_{0}x_{1}-s_{2}s_{0}x_{1})s_{2}y_{2}, \\
C_{(2,0)(1)}(x_{1}\otimes y_{2})&
=(s_{2}s_{0}x_{1}-s_{2}s_{1}x_{1})(s_{1}y_{2}-s_{2}y_{2}) \\
C_{(0)(2,1)}(x_{2}\otimes y_{1})&
=s_{2}s_{1}x_{2}(s_{0}y_{1}-s_{1}y_{1}+s_{2}y_{1}) \\
C_{(1)(0)}(x_{2}\otimes y_{2})&
=[s_{1}x_{2}(s_{0}y_{2}-s_{1}y_{2})+s_{2}(x_{2}y_{2}), \\
C_{(2)(0)}(x_{2}\otimes y_{2})& =(s_{2}x_{2})(s_{0}y_{2}), \\
C_{(2)(1)}(x_{2}\otimes y_{2})& =s_{2}x_{2}(s_{1}y_{2}-s_{2}y_{2}).
\end{align*}

For $n=4$, the key pairings are thus the following

\begin{center}
\begin{tabular}{lllll}
$C_{(3,2,1)(0)},$ & $C_{(3,2,0)(1)},$ & $C_{(3,1,0)(2)},$ & $C_{(2,1,0)(3)},$
& $C_{(3,0)(2,1)},$ \\
$C_{(3,1)(2,0)},$ & $C_{(3,2)(1,0)},$ & $C_{(3,2)(1)},$ & $C_{(3,2)(0)},$ & $%
C_{(3,1)(0)},$ \\
$C_{(0)(2,1)},$ & $C_{(3,1)(2)},$ & $C_{(2,1)(3)},$ & $C_{(3,0)(2)},$ & $%
C_{(3,0)(1)},$ \\
$C_{(2,0)(3)},$ & $C_{(2,0)(1)},$ & $C_{(1,0)(3)},$ & $C_{(1,0)(2)},$ & $%
C_{(3)(2)},$ \\
$C_{(3)(1)},$ & $C_{(3)(0)},$ & $C_{(2)(1)},$ & $C_{(2)(0)},$ & $C_{(1)(0)},$%
\end{tabular}
\end{center}

\begin{theorem}
(\cite{patron3}) Let $\mathbf{E}$ be a simplicial algebra with Moore complex
$\mathbf{NE}$ in which $E_{n}=D_{n},$ is an ideal of $E_{n}$ generated by
the degenerate elements in dimension $n,$ then
\begin{equation*}
\begin{tabular}{l}
$\partial _{n}(NE_{n})=\sum\limits_{I,J}\left[ K_{I},K_{J}\right] $%
\end{tabular}%
\end{equation*}%
\ for $I,J\subseteq \lbrack n-1]$ with $I\cup J=[n-1],$ $I=[n-1]-\{\alpha \}$
$J=[n-1]-\{\beta \}$ where $(\alpha ,\beta )\in P(n)$ for $n=2,3$ and $4,$
\end{theorem}

\bigskip

\begin{remark}
Shortly in \cite{amut4} they defined the normal subgroup $\partial
_{n}(NG_{n}\cap D_{n})$ by $F_{\alpha ,\beta }$ elements which were defined
first by Carrasco in \cite{c}. Castiglioni and Ladra generalised this
inclusion in \cite{ladra}.
\end{remark}

Following \cite{patron3} we have

\begin{lemma}
Let $\mathbf{E}$ be a simplicial algebra with Moore complex $\mathbf{NE}$%
\textbf{\ }of\textbf{\ }length $3$. Then for $n=4$ the images of $C_{\alpha
,\beta }$ elements under $\partial _{4}$ given in Table 1 are trivial.
\end{lemma}

\begin{proof}
Since $NG_{4}=1$ by the results in \cite{patron3} result is trivial.
\end{proof}

\newpage

\begin{tabular}{|l|l|l|l|}
\hline
1 & $d_{4}[C_{(3,2,1)(0)}(x_{1}\otimes y_{3})]$ & $=$ & $%
s_{2}s_{1}x_{1}(s_{0}d_{3}y_{3}-s_{1}d_{3}y_{3}+s_{2}d_{3}y_{3}-y_{3})$ \\
\hline
2 & $d_{4}[C_{(3,2,0)(1)}(x_{1}\otimes y_{3})]$ & $=$ & $%
(s_{2}s_{0}x_{1}-s_{2}s_{1}x_{1})(s_{1}d_{3}y_{3}-s_{2}d_{3}y_{3}+y_{3})$ \\
\hline
3 & $d_{4}[C_{(3,1,0)(2)}(x_{1}\otimes y_{3})]$ & $=$ & $%
(s_{1}s_{0}x_{1}-s_{2}s_{0}x_{1})(s_{2}d_{3}y_{3}-y_{3})$ \\ \hline
4 & $d_{4}[C_{(2,1,0)(3)}(x_{1}\otimes y_{3})]$ & $=$ & $%
(s_{2}s_{1}s_{0}d_{1}x_{1}-s_{1}s_{0}x_{1})y_{3}$ \\ \hline
5 & $d_{4}[C_{(3,2)(1,0)}(x_{2}\otimes y_{2})]$ & $=$ & $%
(s_{1}s_{0}d_{2}x_{2}-s_{2}s_{0}d_{2}x_{2}-s_{0}x_{2})s_{2}y_{2}$ \\ \hline
6 & $d_{4}[C_{(3,1)(2,0)}(x_{2}\otimes y_{2})]$ & $=$ & $%
(s_{1}x_{2}-s_{0}x_{2}+s_{2}s_{0}d_{2}x_{2}-s_{2}s_{1}d_{2}x_{2})(s_{1}y_{2}-s_{2}y_{2})
$ \\ \hline
7 & $d_{4}[C_{(3,0)(2,1)}(x_{2}\otimes y_{2})]$ & $=$ & $%
(s_{2}s_{1}d_{2}x_{2}-s_{1}x_{2})(s_{0}y_{2}-s_{1}y_{2}+s_{2}y_{2})$ \\
\hline
8 & $d_{4}[C_{(3,2)(1)}(x_{2}\otimes y_{3})]$ & $=$ & $%
s_{2}x_{2}(s_{1}d_{3}y_{3}-s_{2}d_{3}y_{3}+y_{3})$ \\ \hline
9 & $d_{4}[C_{(3,2)(0)}(x_{2}\otimes y_{3})]$ & $=$ & $%
s_{2}x_{2}(s_{2}d_{3}y_{3}-s_{1}d_{3}y_{3}+s_{0}d_{3}y_{3}-y_{3})$ \\ \hline
10 & $d_{4}[C_{(3,1)(2)}(x_{2}\otimes y_{3})]$ & $=$ & $%
(s_{1}x_{2}-s_{2}x_{2})(s_{2}d_{3}y_{3}-y_{3})$ \\ \hline
11 & $d_{4}[C_{(3,1)(0)}(x_{2}\otimes y_{3})]$ & $=$ & $%
(s_{1}x_{2}-s_{2}x_{2})(s_{2}d_{3}y_{3}-s_{1}d_{3}y_{3}+s_{0}d_{3}y_{3}-y_{3})
$ \\ \hline
12 & $d_{4}[C_{(3,0)(2)}(x_{2}\otimes y_{3})]$ & $=$ & $%
(s_{0}x_{2}-s_{1}x_{2}+s_{2}x_{2})(s_{2}d_{3}y_{3}-y_{3})$ \\ \hline
13 & $d_{4}[C_{(3,0)(1)}(x_{2}\otimes y_{3})]$ & $=$ & $%
(s_{0}x_{2}-s_{1}x_{2}+s_{2}x_{2})(s_{1}d_{3}y_{3}-s_{2}d_{3}y_{3}+y_{3})$
\\ \hline
14 & $d_{4}[C_{(2,1)(3)}(x_{2}\otimes y_{3})]$ & $=$ & $%
(s_{2}s_{1}d_{2}x_{2}-s_{1}x_{2})y_{3}$ \\ \hline
15 & $d_{4}[C_{(0)(2,1)}(x_{2}\otimes y_{3})]$ & $=$ & $%
(s_{2}s_{1}d_{2}x_{2}-s_{1}x_{2})(s_{2}d_{3}y_{3}-s_{1}d_{3}y_{3}+s_{0}d_{3}y_{3}-y_{3})
$ \\ \hline
16 & $d_{4}[C_{(2,0)(3)}(x_{2}\otimes y_{3})]$ & $=$ & $%
(s_{2}s_{0}d_{2}x_{2}-s_{0}x_{2}+s_{1}x_{2}-s_{1}s_{1}d_{2}x_{2})y_{3}$ \\
\hline
17 & $d_{4}[C_{(2,0)(1)}(x_{2}\otimes y_{3})]$ & $=$ & $%
(s_{2}s_{0}d_{2}x_{2}-s_{0}x_{2}+s_{1}x_{2}-s_{2}s_{1}d_{2}x_{2})$ \\ \hline
&  &  & $(s_{1}d_{3}y_{3}-s_{2}d_{3}y_{3}+y_{3})$ \\ \hline
18 & $d_{4}[C_{(1,0)(3)}(x_{2}\otimes y_{3})]$ & $=$ & $%
(s_{2}s_{0}d_{2}x_{2}-s_{0}x_{2}-s_{1}s_{0}d_{0}x_{2})y_{3}$ \\ \hline
19 & $d_{4}[C_{(1,0)(2)}(x_{2}\otimes y_{3})]$ & $=$ & $%
(s_{1}s_{0}d_{2}x_{2}-s_{2}s_{0}d_{2}x_{2}+s_{0}x_{2})(s_{2}d_{3}y_{3}-y_{3})
$ \\ \hline
20 & $d_{4}[C_{(3)(2)}(x_{3}\otimes y_{3})]$ & $=$ & $%
x_{3}(s_{2}d_{3}y_{3}-y_{3})$ \\ \hline
21 & $d_{4}[C_{(3)(1)}(x_{3}\otimes y_{3})]$ & $=$ & $%
x_{3}(s_{1}d_{3}y_{3}-s_{2}d_{3}y_{3}+y_{3})$ \\ \hline
22 & $d_{4}[C_{(3)(0)}(x_{3}\otimes y_{3})]$ & $=$ & $%
x_{3}(s_{2}d_{3}y_{3}-s_{1}d_{3}y_{3}+s_{0}d_{3}y_{3}-y_{3})$ \\ \hline
23 & $d_{4}[C_{(2)(1)}(x_{3}\otimes y_{3})]$ & $=$ & $%
(s_{2}d_{3}x_{3}-x_{3})(s_{1}d_{3}y_{3}-s_{2}d_{3}y_{3}+y_{3})$ \\ \hline
24 & $d_{4}[C_{(2)(0)}(x_{3}\otimes y_{3})]$ & $=$ & $%
(s_{2}d_{3}x_{3}-x_{3})(s_{2}d_{3}y_{3}-s_{1}d_{3}y_{3}+s_{0}d_{3}y_{3}-y_{3})
$ \\ \hline
25 & $d_{4}[C_{(1)(0)}(x_{3}\otimes y_{3})]$ & $=$ & $%
(s_{1}d_{3}x_{3}-s_{2}d_{3}x_{3}+x_{3})(s_{2}d_{3}y_{3}-s_{1}d_{3}y_{3}+s_{0}d_{3}y_{3}-y_{3})
$ \\ \hline
\end{tabular}

\begin{center}
Table 1
\end{center}

where $x_{3},y_{3}\in NG_{3},x_{2},y_{2}\in NG_{2},x_{1}\in NG_{1}$.

\newpage

\subsection{Crossed Modules}

Here we will recall the notion of crossed modules of commutative algebras
given in \cite{portervascon} and \cite{e1}

Let $R$ be a $\mathbf{k}$-algebra with identity. A $crossed$ $module$ of
commutative algebra is an $R$-algebra $C$, together with a commutative
action of $R$ on $C$ and $R$-algebra morphism $\partial :C\rightarrow R$
together with an action of $R$ on $C$, written $r\cdot c$ for $r\in R$ and $%
c\in C$, satisfying the conditions.

\textbf{CM1)} for all $r\in R$ , $c\in C $%
\begin{equation*}
\partial (r\cdot c)=r\partial c
\end{equation*}

\textbf{CM2) }(Peiffer Identity) for all $c,c^{\prime }\in C$%
\begin{equation*}
{\partial c}\cdot c^{\prime }=cc^{\prime }
\end{equation*}%
We will denote such a crossed module by $(C,R,\partial )$.

A \textit{morphism of crossed module} from $(C,R,\partial )$ to $(C^{\prime
},R^{\prime },\partial ^{\prime })$ is a pair of $\mathbf{k}$-algebra
morphisms%
\begin{equation*}
\phi :C\longrightarrow C^{\prime }\text{ , \ }\psi :R\longrightarrow
R^{\prime }\text{ }
\end{equation*}%
such that $\phi (r\cdot c)={\psi (r)}\cdot \phi (r)$ and $\partial ^{\prime
}\phi (c)=\psi \partial (c)$.

We thus get a category $\mathbf{XMod}$ of crossed modules. \newline

\textit{Examples of Crossed Modules}

(i) Any ideal, $I$, in $R$ gives an inclusion map
$I\longrightarrow R,$ which is a crossed module then we will say
$\left( I,R,i\right) $ is \ an ideal pair. In this case, of
course, $R$ acts on $I$ by multiplication and the inclusion
homomorphism $i$ makes $\left( I,R,i\right) $ into a crossed
module, an \textquotedblleft inclusion crossed
modules\textquotedblright . Conversely,

\begin{lemma}
If $(C,R,\partial )$ is a crossed module, $\partial (C)$ is an
ideal of $R.$
\end{lemma}

(ii) Any $R$-module $M$ can be considered as an $R$-algebra with
zero multiplication and hence the zero morphism $0:M\rightarrow R$
sending everything in $M$ to the zero element of $R$ is a crossed
module. Again conversely:

\begin{lemma}
If $(C,R,\partial )$ is a crossed module, $\ker \partial $ is an
ideal in $C$ and inherits a natural $R$-module structure from
$R$-action on $C.$ Moreover, $\partial (C)$ acts trivially on
$\ker \partial ,$ hence $\ker
\partial $ has a natural $R/\partial (C)$-module structure.
\end{lemma}

As these two examples suggest, general crossed modules lie between
the two extremes of ideal and modules. Both aspects are
important.\bigskip

(iii) In the category of algebras, the appropriate replacement for
automorphism \ groups is the multiplication algebra defined by Mac
Lane \cite{[m]}. Then automorphism crossed module correspond to
the multiplication crossed module $\left( R,M\left( R\right) ,\mu
\right) $.

To see this crossed module, we need to assume $Ann\left( R\right) =0$ or $%
R^{2}=R$ and let $M\left( R\right) $ be the set of all multipliers
$\delta :R\rightarrow R$ such that for all $c,c^{\prime }\in C$,\
$\delta \left( rr^{\prime }\right) =\delta \left( r\right)
r^{\prime }.$ $M\left( R\right) $ acts on $R$ by
\[
\begin{array}{rcl}
M\left( R\right) \times R & \longrightarrow & R \\
\left( \delta ,r\right) & \longmapsto & \delta \left( r\right)
\end{array}
\]
and there is a morphism $\mu :R\rightarrow M\left( R\right) $
defined by $\mu \left( r\right) =\delta _{r}$ with $\delta
_{r}\left( r^{\prime }\right) =rr^{\prime }$ for all $r,r^{\prime
}\in R.$

\subsection{$2$-Crossed Modules}

Now we recall the commutative algebra case of $2$-crossed modules due to
A.R.Grandjean and Vale, \cite{grand}.

A 2-crossed module of $k$-algebras is a complex
\begin{equation*}
C_{2}\overset{\partial _{2}}{\longrightarrow }C_{1}\overset{\partial _{1}}{%
\longrightarrow }C_{0}
\end{equation*}%
of $C_{0}$-algebras with $\partial _{2},\partial _{1}$ morphisms of $C_{0}$%
-algebras, where $C_{0}$ acts on $C_{0}$ by multiplication, with a bilinear
function
\begin{equation*}
\{\quad \otimes \quad \}:C_{1}\otimes _{C_{0}}C_{1}\longrightarrow C_{2}
\end{equation*}%
called as \textit{Peiffer lifting} which satisfies the following axioms:
\begin{equation*}
\begin{array}{lrrl}
\mathbf{2CM1)} & \partial _{2}\{y_{0}\otimes y_{1}\} & = & y_{0}y_{1}-\text{
}^{\partial _{1}(y_{1})}y_{0} \\
\mathbf{2CM2)} & \{\partial _{2}(x_{1})\otimes \partial _{2}(x_{2})\} & = &
x_{1}x_{2} \\
\mathbf{2CM3)} & \{y_{0}\otimes y_{1}y_{2}\} & = & \{y_{0}y_{1}\otimes
y_{2}\}+\text{ }^{\partial _{1}y_{2}}\{{y_{0}\otimes y_{1}}\} \\
\mathbf{2CM4)} & (i)\quad \{\partial _{2}(x)\otimes y\} & = & y\cdot x-\text{
}^{\partial _{1}(y)}x \\
& (ii)\quad \{y\otimes {\partial _{2}(x)}\} & = & y\cdot x\newline
\\
\mathbf{2CM5)} & \text{ }^{z}\{y_{0}\otimes y_{1}\} & = & \{^{z}y_{0}\otimes
y_{1}\}=\{y_{0}\text{ }\otimes \text{ }^{z}y_{1}\}\newline
\end{array}%
\end{equation*}%
\newline
for all $x,x_{1},x_{2}\in C_{2}$, $y,y_{0},y_{1},y_{2}\in C_{1}$ and $z\in
C_{0}$. A morphism of $2$-crossed modules can be defined in an obvious way.
We thus define the category of $2$-crossed modules denoting it by $\mathbf{X}%
_{2}\mathbf{Mod}$.

The proof of the following theorem can be found in \cite{patron3}.

\begin{theorem}
The category of $2$-crossed modules is equivalent to the category of
simplicial algebras with Moore complex of length $2$.
\end{theorem}

Now we will give some remarks on $2$-crossed modules where the group case
can be found in \cite{timcrossed}.

1) Let $C_{1}\overset{\partial _{1}}{\longrightarrow }C_{0}$ be a crossed
module. If we take $C_{2}$ trivial then
\begin{equation*}
C_{2}\overset{\partial _{2}}{\longrightarrow }C_{1}\overset{\partial _{1}}{%
\longrightarrow }C_{0}
\end{equation*}%
is a $2$-crossed module with the Peiffer lifting defined by $\left\{
x\otimes y\right\} =0$ for $x,y\in C_{1}$.

2) If
\begin{equation*}
C_{2}\overset{\partial _{2}}{\longrightarrow }C_{1}\overset{\partial _{1}}{%
\longrightarrow }C_{0}
\end{equation*}%
is a $2$-crossed module then
\begin{equation*}
\frac{C_{1}}{{Im}\partial _{2}}\overset{\partial _{1}}{\longrightarrow }%
C_{0}
\end{equation*}%
is a crossed module.

3) Let
\begin{equation*}
C_{2}\overset{\partial _{2}}{\longrightarrow }C_{1}\overset{\partial _{1}}{%
\longrightarrow }C_{0}
\end{equation*}%
be a $2$-crossed module with trivial Peiffer lifting then $C_{1}\overset{%
\partial _{1}}{\longrightarrow }C_{0}$ will be a crossed module. Also in
this situation we have the trivial action of $C_{0}$ on $C_{2}$.

\section{Three Crossed Modules}

\qquad As a consequence of \cite{arkaus}, here we will define $3$-crossed
modules of commutative algebras. The way is similar but some of the
conditions are different.

Let $\mathbf{E}$ be a simplicial algebra with Moore complex of length $3$
and $NE_{0}=C_{0},$ $NE_{1}=C_{1},$ $NE_{2}=C_{2},$ $NE_{3}=C_{3}$. Thus we
have a $\mathbf{k}$-algebra complex
\begin{equation*}
C_{3}\overset{\partial _{3}}{\longrightarrow }C_{2}\overset{\partial _{2}}{%
\longrightarrow }C_{1}\overset{\partial _{1}}{\longrightarrow }C_{0}
\end{equation*}%
Let the actions of $C_{0}$ on $C_{3}$, $C_{2}$, $C_{1}$, $C_{1}$ on $C_{2}$,
$C_{3}$ and $C_{2}$ on $C_{3}$ be as follows;%
\begin{equation*}
\begin{array}{ll}
^{x_{0}}x_{1} & =s_{0}x_{0}x_{1} \\
^{x_{0}}x_{2} & =s_{1}s_{0}x_{0}x_{2} \\
^{x_{0}}x_{3} & =s_{2}s_{1}s_{0}x_{0}x_{3} \\
^{x_{1}}x_{2} & =s_{1}x_{1}x_{2} \\
^{x_{1}}x_{3} & =s_{2}s_{1}x_{1}x_{3} \\
x_{2}\cdot x_{3} & =s_{2}x_{2}x_{3}%
\end{array}%
\end{equation*}

\bigskip

Then, since%
\begin{equation*}
\begin{array}{rr}
(s_{2}s_{1}s_{0}\partial _{1}x_{1}-s_{1}s_{0}x_{1})y_{3} & =0 \\
(s_{2}s_{1}\partial _{2}x_{2}-s_{1}x_{2})y_{3} & =0 \\
x_{3}(s_{2}\partial _{3}y_{3}-y_{3}) & =0%
\end{array}%
\end{equation*}%
we get

\begin{equation*}
\begin{array}[t]{ll}
^{\partial _{1}x_{1}}\text{ }y_{3} & =s_{1}s_{0}x_{1}y_{3} \\
^{\partial _{2}x_{2}}\text{ }y_{3} & =s_{1}x_{2}y_{3} \\
\partial _{3}x_{3}\cdot y_{3} & =x_{3}y_{3}%
\end{array}%
\end{equation*}%
and using the simplicial identities we get,
\begin{equation*}
\partial _{3}(x_{2}\cdot x_{3})=\partial _{3}(s_{2}x_{2}x_{3})=\partial
_{3}(s_{2}x_{2})\partial _{3}(x_{3})=x_{2}\partial _{3}(x_{3})
\end{equation*}%
Thus $\partial _{3}:C_{3}\rightarrow C_{2}$ is a crossed module.

\begin{definition}
Let $C_{3}\overset{\partial _{3}}{\longrightarrow }C_{2}\overset{\partial
_{2}}{\longrightarrow }C_{1}\overset{\partial _{1}}{\longrightarrow }C_{0}$
be a complex of $\mathbf{k}$-algebras defined above. We define Peiffer
liftings as follows;
\begin{equation*}
\begin{array}{lllll}
\left\{ ~\otimes ~\right\} & : & C_{1}\otimes C_{1} & \rightarrow & C_{2} \\
&  & \left\{ x_{1}\otimes y_{1}\right\} & = &
s_{1}x_{1}(s_{0}y_{1}-s_{1}y_{1}) \\
\left\{ ~\otimes ~\right\} _{(1,0)(2)} & : & C_{1}\otimes C_{2} & \rightarrow
& C_{3} \\
&  & \left\{ x_{1}\otimes y_{2}\right\} _{(1,0)(2)} & = &
(s_{1}s_{0}x_{1}-s_{2}s_{0}x_{1})s_{2}y_{2} \\
\left\{ ~\otimes ~\right\} _{(2,0)(1)} & : & C_{1}\otimes C_{2} & \rightarrow
& C_{3} \\
&  & \left\{ x_{1}\otimes y_{2}\right\} _{(2,0)(1)} & = &
(s_{2}s_{0}x_{1}-s_{2}s_{1}x_{1})(s_{1}y_{2}-s_{2}y_{2}) \\
\left\{ ~\otimes ~\right\} _{(0)(2,1)} & : & C_{1}\otimes C_{2} & \rightarrow
& C_{3} \\
&  & \left\{ x_{1}\otimes y_{2}\right\} _{(0)(2,1)} & = &
s_{2}s_{1}x_{1}(s_{0}y_{2}-s_{1}y_{2}+s_{2}y_{2}) \\
\left\{ ~\otimes ~\right\} _{(1)(0)} & : & C_{2}\otimes C_{2} & \rightarrow
& C_{3} \\
&  & \left\{ x_{2}\otimes y_{2}\right\} _{(1)(0)} & = &
(s_{0}x_{2}-s_{1}x_{2})s_{1}y_{2}+s_{2}(x_{2}y_{2}) \\
\left\{ ~\otimes ~\right\} _{(2)(0)} & : & C_{2}\otimes C_{2} & \rightarrow
& C_{3} \\
&  & \left\{ x_{2}\otimes y_{2}\right\} _{(2)(0)} & = & s_{2}x_{2}s_{0}y_{2}
\\
\left\{ ~\otimes ~\right\} _{(2)(1)} & : & C_{2}\otimes C_{2} & \rightarrow
& C_{3} \\
&  & \left\{ x_{2}\otimes y_{2}\right\} _{(2)(1)} & = &
s_{2}x_{2}(s_{1}y_{2}-s_{2}y_{2})%
\end{array}%
\end{equation*}%
where \ $x_{1}$, $y_{1}\in C_{1},$ $x_{2}$ $y_{2}\in C_{2}$.
\end{definition}

Then using Table 1 we get the following identities.

\begin{center}
\begin{tabular}{|l|l|l|}
\hline
$\left\{ x_{2}\otimes \partial _{2}y_{2}\right\} _{(0)(2,1)}$ & $=$ & $%
\left\{ x_{2}\otimes y_{2}\right\} _{(1)(0)}+\left\{ x_{2}\otimes
y_{2}\right\} _{(2)(1)}$ \\ \hline
$\left\{ x_{1}\otimes \partial _{3}y_{3}\right\} _{(2,0)(1)}$ & $=$ & $%
\left\{ x_{1}\otimes \partial _{3}y_{3}\right\} _{(0)(2,1)}{\small +}\left\{
x_{1}\otimes \partial _{3}y_{3}\right\} _{(1,0)(2)}{\small -}$ $^{{\small %
\partial }_{1}{\small x}_{1}}{\small y}_{3}$ \\ \hline
$\left\{ \partial _{2}x_{2}\otimes y_{2}\right\} _{(1,0)(2)}$ & $=$ & $%
-\left\{ x_{2}\otimes y_{2}\right\} _{(0)(2)}$ \\ \hline
$\left\{ \partial _{3}x_{3}\otimes \partial _{3}y_{3}\right\} _{(1)(0)}$ & $%
= $ & $x_{3}y_{3}$ \\ \hline
$\left\{ \partial _{2}x_{2}\otimes \partial _{3}y_{3}\right\} _{(0)(2,1)}$ &
$=$ & $^{\partial _{2}x_{2}}y_{3}$ \\ \hline
$\left\{ \partial _{2}x_{2}\otimes \partial _{3}y_{3}\right\} _{(1,0)(2)}$ &
$=$ & $-\left\{ x_{2}\otimes \partial _{3}y_{3}\right\} _{(0)(2)}$ \\ \hline
$\left\{ \partial _{2}x_{2}\otimes \partial _{3}y_{3}\right\} _{(2,0)(1)}$ &
$=$ & $\partial _{2}x_{2}\cdot y_{3}-\left\{ x_{2}\otimes \partial
_{3}y_{3}\right\} _{(0)(2)}$ \\ \hline
$\left\{ x_{2}\otimes y_{2}y\prime _{2}\right\} _{(1)(0)}$ & $=$ & $\left\{
x_{2}\otimes \partial _{2}(y_{2}y\prime _{2})\right\} _{(0)(2,1)}{\small -}%
\left\{ x_{2}\otimes (y_{2}y\prime _{2})\right\} _{(2)(1)}$ \\ \hline
$\left\{ x\prime _{2}x_{2}\otimes y_{2}\right\} _{(1)(0)}$ & $=$ & $\left\{
x\prime _{2}x_{2}\otimes \partial _{2}y_{2}\right\} _{(0)(2,1)}-\left\{
x\prime _{2}x_{2}\otimes y_{2}\right\} _{(2)(1)}$ \\ \hline
$\left\{ x_{2}\otimes y_{2}y\prime _{2}\right\} _{(2)(1)}$ & $=$ & $\left\{
x_{2}\otimes \partial _{2}(y_{2}y\prime _{2})\right\} _{(1)(2,0)}+\left\{
x_{2}\otimes y_{2}y\prime _{2}\right\} _{(2)(0)}-\left\{ x_{2}\otimes
y_{2}y\prime _{2}\right\} _{(1)(0)}$ \\ \hline
$\left\{ x_{2}x\prime _{2}\otimes y_{2}\right\} _{(2)(1)}$ & $=$ & $\left\{
x_{2}x\prime _{2}\otimes \partial _{2}y_{2}\right\} _{(1)(2,0)}+\left\{
x_{2}x\prime _{2}\otimes y_{2}\right\} _{(2)(0)}-\left\{ x_{2}x\prime
_{2}\otimes y_{2}\right\} _{(1)(0)}$ \\ \hline
$\left\{ x_{2}\otimes y_{2}y\prime _{2}\right\} _{(2)(0)}$ & $=$ & $-\left\{
x_{2}\otimes \partial _{2}(y_{2}y\prime _{2})\right\} _{(2)(1,0)}$ \\ \hline
$\left\{ x_{2}\otimes \partial _{3}y_{3}\right\} _{(2)(1)}$ & $=$ & $%
x_{2}\cdot y_{3}$ \\ \hline
$\left\{ \partial _{3}x_{3}\otimes y_{2}\right\} _{(2)(1)}$ & $=$ & $%
x_{3}^{\partial _{2}y_{2}}+x_{3}\cdot y_{2}$ \\ \hline
$\left\{ \partial _{3}x_{3}\otimes y_{2}\right\} _{(1)(0)}$ & $=$ & $%
y_{2}\cdot x_{3}$ \\ \hline
$\left\{ \partial _{3}x_{3}\otimes y_{2}\right\} _{(2)(0)}$ & $=$ & $0$ \\
\hline
$\partial _{3}\left\{ x_{2}\otimes y_{2}\right\} _{(2)(0)}$ & $=$ & $%
-\partial _{3}\left\{ \partial _{2}x_{2}\otimes y_{2}\right\} _{(1,0)(2)}$
\\ \hline
$\partial _{3}\left\{ x_{2}\otimes y_{2}\right\} _{(1)(0)}$ & $=$ & $\left\{
\partial _{2}x_{2}\otimes \partial _{2}y_{2}\right\} _{(1)(0)}+x_{2}y_{2}$
\\ \hline
$\partial _{3}\left\{ x_{2}\otimes y_{2}\right\} _{(2)(1)}$ & $=$ & $%
^{\partial _{2}y_{2}}x_{2}-x_{2}y_{2}$ \\ \hline
$\partial _{3}\left\{ x_{1}\otimes y_{2}\right\} _{(2,0)(1)}$ & $=$ & $%
\partial _{3}\left\{ x_{1}\otimes y_{2}\right\} _{(1,0)(2)}+\left\{
x_{1}\otimes \partial _{2}y_{2}\right\} -^{\partial
_{1}x_{1}}y_{2}+^{x_{1}}y_{2}$ \\ \hline
$\partial _{3}\left\{ x_{1}\otimes y_{2}\right\} _{(0)(2,1)}$ & $=$ & $%
\left\{ x_{1}\otimes \partial _{2}y_{2}\right\} +^{x_{1}}y_{2}$ \\ \hline
\end{tabular}

Table 2

\begin{tabular}{|l|l|l|l|l|}
\hline
$^{x_{0}}\left\{ x_{1}\otimes y_{1}\right\} $ & $=$ & $\left\{
^{x_{0}}x_{1}\otimes y_{1}\right\} $ & $=$ & $\left\{ x_{1}\otimes
^{x_{0}}y_{1}\right\} $ \\ \hline
$^{x_{0}}\left\{ x_{1}\otimes y_{2}\right\} _{(1,0)(2)}$ & $=$ & $\left\{
^{x_{0}}x_{1}\otimes y_{2}\right\} _{(1,0)(2)}$ & $=$ & $\left\{
x_{1}\otimes ^{x_{0}}y_{2}\right\} _{(1,0)(2)}$ \\ \hline
$^{x_{0}}\left\{ x_{1}\otimes y_{2}\right\} _{(0)(2,1)}$ & $=$ & $\left\{
^{x_{0}}x_{1}\otimes y_{2}\right\} _{(0)(2,1)}$ & $=$ & $\left\{
x_{1}\otimes ^{x_{0}}y_{2}\right\} _{(0)(2,1)}$ \\ \hline
$^{x_{0}}\left\{ x_{1}\otimes y_{2}\right\} _{(2,0)(1)}$ & $=$ & $\left\{
^{x_{0}}x_{1}\otimes y_{2}\right\} _{(2,0)(1)}$ & $=$ & $\left\{
x_{1}\otimes ^{x_{0}}y_{2}\right\} _{(2,0)(1)}$ \\ \hline
$^{x_{0}}\left\{ x_{2}\otimes y_{2}\right\} _{(1)(0)}$ & $=$ & $\left\{
^{x_{0}}x_{2}\otimes y_{2}\right\} _{(1)(0)}$ & $=$ & $\left\{ x_{2}\otimes
^{x_{0}}y_{2}\right\} _{(1)(0)}$ \\ \hline
$^{x_{0}}\left\{ x_{2}\otimes y_{2}\right\} _{(2)(0)}$ & $=$ & $\left\{
^{x_{0}}x_{2}\otimes y_{2}\right\} _{(2)(0)}$ & $=$ & $\left\{ x_{2}\otimes
^{x_{0}}y_{2}\right\} _{(2)(0)}$ \\ \hline
$^{x_{0}}\left\{ x_{2}\otimes y_{2}\right\} _{(2)(1)}$ & $=$ & $\left\{
^{x_{0}}x_{2}\otimes y_{2}\right\} _{(2)(1)}$ & $=$ & $\left\{ x_{2}\otimes
^{x_{0}}y_{2}\right\} _{(2)(1)}$ \\ \hline
\end{tabular}

Table 3

\begin{tabular}{|l|l|l|l|l|}
\hline
$^{z_{1}}\left\{ x_{1}\otimes y_{1}\right\} $ & $=$ & $\left\{
^{z_{1}}x_{1}\otimes y_{1}\right\} $ & $=$ & $\left\{ x_{1}\otimes
^{z_{1}}y_{1}\right\} $ \\
$^{z_{1}}\left\{ x_{1}\otimes y_{2}\right\} _{(1,0)(2)}$ & $=$ & $\left\{
^{z_{1}}x_{1}\otimes y_{2}\right\} _{(1,0)(2)}$ & $=$ & $\left\{
x_{1}\otimes ^{z_{1}}y_{2}\right\} _{(1,0)(2)}$ \\
$^{z_{1}}\left\{ x_{1}\otimes y_{2}\right\} _{(0)(2,1)}$ & $=$ & $\left\{
^{z_{1}}x_{1}\otimes y_{2}\right\} _{(0)(2,1)}$ & $=$ & $\left\{
x_{1}\otimes ^{z_{1}}y_{2}\right\} _{(0)(2,1)}$ \\
$^{z_{1}}\left\{ x_{1}\otimes y_{2}\right\} _{(2,0)(1)}$ & $=$ & $\left\{
^{z_{1}}x_{1}\otimes y_{2}\right\} _{(2,0)(1)}$ & $=$ & $\left\{
x_{1}\otimes ^{z_{1}}y_{2}\right\} _{(2,0)(1)}$ \\
$^{z_{1}}\left\{ x_{2}\otimes y_{2}\right\} _{(1)(0)}$ & $=$ & $\left\{
^{z_{1}}x_{2}\otimes y_{2}\right\} _{(1)(0)}$ & $=$ & $\left\{ x_{2}\otimes
^{z_{1}}y_{2}\right\} _{(1)(0)}$ \\
$^{z_{1}}\left\{ x_{2}\otimes y_{2}\right\} _{(2)(0)}$ & $=$ & $\left\{
^{z_{1}}x_{2}\otimes y_{2}\right\} _{(2)(0)}$ & $=$ & $\left\{ x_{2}\otimes
^{z_{1}}y_{2}\right\} _{(2)(0)}$ \\
$^{z_{1}}\left\{ x_{2}\otimes y_{2}\right\} _{(2)(1)}$ & $=$ & $\left\{
^{z_{1}}x_{2}\otimes y_{2}\right\} _{(2)(1)}$ & $=$ & $\left\{ x_{2}\otimes
^{z_{1}}y_{2}\right\} _{(2)(1)}$ \\ \hline
\end{tabular}

Table 4
\end{center}

\noindent where $x_{0}\in C_{0}$, $x_{1},y_{1}\in C_{1},$ $x_{2},y_{2}\in C_{2},$ $%
x_{3},y_{3}\in C_{3}$. From these results all liftings given in definition 1
are $C_{0}$,$C_{1}$-bilinear maps.

\begin{definition}
A \textit{3-crossed module} consist of a complex
\begin{equation*}
C_{3}\overset{\partial _{3}}{\longrightarrow }C_{2}\overset{\partial _{2}}{%
\longrightarrow }C_{1}\overset{\partial _{1}}{\longrightarrow }C_{0}
\end{equation*}%
together with $\partial _{3}$, $\partial _{2}$,$\partial _{1}$ which are $%
C_{0},C_{1}$-algebra morphisms, an action of $C_{0}$ on $C_{3},C_{2},C_{1}$,
an action of $C_{1}$ on $C_{2},C_{3}$ and an action of $C_{2}$ on $C_{3},$%
further $C_{0},C_{1}$-bilinear maps
\begin{equation*}
\begin{tabular}{lll}
$\{$ $\otimes $ $\}_{(1)(0)}:C_{2}\otimes C_{2}\longrightarrow C_{3},$ & $\{$
$\otimes $ $\}_{(0)(2)}:C_{2}\otimes C_{2}\longrightarrow C_{3},$ & $\{$ $%
\otimes $ $\}_{(2)(1)}:C_{2}\otimes C_{2}\longrightarrow C_{3},$ \\
&  &  \\
$\{$ $\otimes $ $\}_{(1,0)(2)}:C_{1}\otimes C_{2}\longrightarrow C_{3},$ & $%
\{$ $\otimes $ $\}_{(2,0)(1)}:C_{1}\otimes C_{2}\longrightarrow C_{3},$ &
\\
&  &  \\
$\{$ $\otimes $ $\}_{(0)(2,1)}:C_{2}\otimes C_{1}\longrightarrow C_{3},$ & $%
\{$ $\otimes $ $\}:C_{1}\otimes C_{1}\longrightarrow C_{2}$ &
\end{tabular}%
\end{equation*}%
called \textit{\ Peiffer lifting}s which satisfy the following axioms for
all $x_{1}\in C_{1},$ $x_{2},y_{2}\in C_{2},$ and $x_{3},y_{3}\in C_{3}$:

\begin{equation*}
\begin{array}{lrrl}
\mathbf{3CM1)} &  & C_{3}\overset{\partial _{3}}{\longrightarrow }C_{2}%
\overset{\partial _{2}}{\longrightarrow }C_{1} & \text{is a }2\text{-crossed
module with the Peiffer lifting }\{\text{ }\otimes \text{ }\}_{(2),(1)} \\
\mathbf{3CM2)} &  & \partial _{2}\left\{ x_{1}\otimes y_{1}\right\} & =\text{
}^{\text{ }\partial _{1}y_{1}}x_{1}-x_{1}y_{1} \\
\mathbf{3CM3)} &  & \left\{ x_{2}\otimes \partial _{2}y_{2}\right\}
_{(0)(2,1)} & =\left\{ x_{2}\otimes y_{2}\right\} _{(2)(1)}-\left\{
x_{2}\otimes y_{2}\right\} _{(1)(0)} \\
\mathbf{3CM4)} &  & \partial _{3}\left\{ x_{2}\otimes y_{2}\right\} _{(1)(0)}
& =\left\{ \partial _{2}x_{2}\otimes \partial _{2}y_{2}\right\} +x_{2}y_{2}
\\
\mathbf{3CM5)} &  & \left\{ x_{1}\otimes \partial _{3}y_{3}\right\}
_{(2,0)(1)} & =\left\{ x_{1}\otimes \partial _{3}y_{3}\right\}
_{(0)(2,1)}+\left\{ x_{1}\otimes \partial _{3}y_{3}\right\} _{(1,0)(2)}-%
\text{ }^{\partial _{1}x_{1}}y_{3} \\
\mathbf{3CM6)} &  & \left\{ \partial _{2}x_{2}\otimes y_{2}\right\}
_{(2,0)(1)} & =-\left\{ x_{2}\otimes y_{2}\right\} _{(0)(2)}+\left(
x_{2}y_{2}\right) \cdot \left\{ x_{2}\otimes y_{2}\right\} _{(2)(1)}+\left\{
x_{2}\otimes y_{2}\right\} _{(1)(0)} \\
\mathbf{3CM7)} &  & \left\{ \partial _{3}x_{3}\otimes \partial
_{3}y_{3}\right\} _{(1)(0)} & =y_{3}x_{3} \\
\mathbf{3CM8}) &  & \left\{ \partial _{3}y_{3}\otimes \partial
_{2}x_{2}\right\} _{(0)(2,1)} & =\text{ }^{-\partial _{2}x_{2}}y_{3} \\
\mathbf{3CM9)} &  & \left\{ \partial _{2}x_{2}\otimes \partial
_{3}y_{3}\right\} _{(1,0)(2)} & =-\left\{ x_{2}\otimes \partial
_{3}y_{3}\right\} _{(0)(2)} \\
\mathbf{3CM10)} &  & \left\{ \partial _{2}x_{2}\otimes \partial
_{3}y_{3}\right\} _{(2,0)(1)} & =\text{ }^{\partial _{2}x_{2}}y_{3}-\left\{
x_{2}\otimes \partial _{3}y_{3}\right\} _{(0)(2)} \\
\mathbf{3CM11)} &  & \left\{ \partial _{3}y_{3}\otimes x_{1}\right\}
_{(0)(2,1)} & =\text{ }^{-x_{1}}y_{3} \\
\mathbf{3CM12)} &  & \left\{ y_{2}\otimes \partial _{3}x_{3}\right\}
_{(1)(0)} & =-y_{2}\cdot x_{3} \\
\mathbf{3CM13)} &  & \left\{ \partial _{3}x_{3}\otimes y_{2}\right\}
_{(1)(0)} & =y_{2}\cdot x_{3} \\
\mathbf{3CM14)} &  & \left\{ \partial _{3}x_{3}\otimes y_{2}\right\}
_{(2)(0)} & =0 \\
\mathbf{3CM15)} &  & \partial _{3}\left\{ x_{1}\otimes y_{2}\right\}
_{(2,0)(1)} & =\partial _{3}\left\{ x_{1}\otimes y_{2}\right\}
_{(1,0)(2)}+\left\{ x_{1}\otimes \partial _{2}y_{2}\right\} -\text{ }%
^{\partial _{1}x_{1}}y_{2}+\text{ }^{x_{1}}y_{2} \\
\mathbf{3CM16)} &  & \partial _{3}\left\{ x_{1}\otimes y_{2}\right\}
_{(0)(2,1)} & =\left\{ x_{1}\otimes \partial _{2}y_{2}\right\} -\text{ }%
^{x_{1}}y_{2}%
\end{array}%
\end{equation*}
\end{definition}

We denote such a 3-crossed module by $(C_{3},C_{2},C_{1},C_{0},\partial
_{3},\partial _{2},\partial _{1}).$

A \textit{morphism of $3$-crossed modules} of groups may be pictured by the
diagram

\begin{center}
$\xymatrix@R=40pt@C=40pt{\ C_{3} \ar[d]_-{f_{3}} \ar@{->}@<0pt>[r]^-{%
\partial_{3}} & C_{2} \ar[d]_-{f_{2}} \ar@{->}@<0pt>[r]^-{\partial_{2}} &
C_{1} \ar[d]_-{f_{1}} \ar@{->}@<0pt>[r]^-{\partial_{1}} & C_{0} \ar[d]_-{%
f_{0}} \\
C^{\prime }_{3} \ar@{->}@<0pt>[r]_-{\partial^{\prime }_{3}} & C^{\prime
}_{2} \ar@{->}@<0pt>[r]_-{\partial^{\prime }_{2}} & C^{\prime }_{1} \ar@{->}%
@<0pt>[r]_-{\partial^{\prime }_{1}} & C^{\prime }_{0} } $ 
\end{center}

\noindent where
\begin{equation*}
f_{1}(^{c_{0}}c_{1})=\text{ }^{(f_{0}(c_{0}))}f_{1}(c_{1}),\text{ }%
f_{2}(^{c_{0}}c_{2})=\text{ }^{(f_{0}(c_{0}))}f_{2}(c_{2}),\text{ }%
f_{3}(^{c_{0}}c_{3})=\text{ }^{(f_{0}(c_{0}))}f_{3}(c_{3})
\end{equation*}%
for $\left\{ \text{ }\otimes \text{ }\right\} _{(0)(2)},\left\{ \text{ }%
\otimes \text{ }\right\} _{(2)(1)},$ $\left\{ \text{ }\otimes \text{ }%
\right\} _{(1)(0)}$
\begin{equation*}
\left\{ \text{ }\otimes \text{ }\right\} f_{2}\otimes f_{2}=f_{3}\left\{
\text{ }\otimes \text{ }\right\}
\end{equation*}%
for $\left\{ \text{ }\otimes \text{ }\right\} _{(1,0)(2)},\left\{ \text{ }%
\otimes \text{ }\right\} _{(2,0)(1)}$ ,$\left\{ \text{ }\otimes \text{ }%
\right\} _{(0)(2,1)}$
\begin{equation*}
\left\{ \text{ }\otimes \text{ }\right\} f_{1}\otimes f_{2}=f_{3}\left\{
\text{ }\otimes \text{ }\right\}
\end{equation*}%
for $\left\{ \text{ }\otimes \text{ }\right\} $
\begin{equation*}
\left\{ \text{ }\otimes \text{ }\right\} f_{1}\otimes f_{1}=f_{2}\left\{
\text{ }\otimes \text{ }\right\}
\end{equation*}%
for all $c_{3}\in C_{3},c_{2}\in C_{3},c_{1}\in C_{3},c_{0}\in C_{3}$. These
compose in an obvious way. So we can define the category of $3$-crossed
modules of commutative algebras,which we will be denoted by $\mathbf{X}_{3}%
\mathbf{ModAlg}$.

\section{Applications}

\subsection{Simplicial Algebras}

As an application we consider the relation between simplicial algebras and $%
3 $-crossed modules which were given in \cite{arkaus} for group case. So
proofs in this section are omitted, since can be checked easily by using the
proofs given in \cite{arkaus}.

\begin{proposition}
Let $\mathbf{E}$ be a simplicial algebra with Moore complex $\mathbf{NE}$.
Then the complex
\begin{equation*}
NE_{3}/\partial _{4}(NE_{4}\cap D_{4})\overset{\overline{\partial }_{3}}{%
\longrightarrow }NE_{2}\overset{\partial _{2}}{\longrightarrow }NE_{1}%
\overset{\partial _{1}}{\longrightarrow }NE_{0}
\end{equation*}%
is a $3$-crossed module with the Peiffer liftings defined below:
\begin{equation*}
\begin{array}{lllll}
\left\{ ~\otimes ~\right\} & : & NE_{1}\otimes NE_{1} & \longrightarrow &
NE_{2} \\
&  & \left\{ x_{1}\otimes y_{1}\right\} _{(1)(0)} & \longmapsto & \overline{%
s_{1}x_{1}(s_{1}y_{1}-s_{0}y_{1})} \\
\left\{ ~\otimes ~\right\} _{(1,0)(2)} & : & NE_{1}\otimes NE_{2} &
\longrightarrow & NE_{3}/\partial _{4}(NE_{4}\cap D_{4}) \\
&  & \left\{ x_{1}\otimes y_{2}\right\} _{(1,0)(2)} & \longmapsto &
\overline{(s_{2}s_{0}x_{1}-s_{1}s_{0}x_{1})s_{2}y_{2}} \\
\left\{ ~\otimes ~\right\} _{(2,0)(1)} & : & NE_{1}\otimes NE_{2} &
\longrightarrow & NE_{3}/\partial _{4}(NE_{4}\cap D_{4})C_{3} \\
&  & \left\{ x_{1}\otimes y_{2}\right\} _{(2,0)(1)} & \longmapsto &
\overline{(s_{2}s_{1}x_{1}-s_{2}s_{0}x_{1})(s_{1}y_{2}-s_{2}y_{2})} \\
\left\{ ~\otimes ~\right\} _{(0)(2,1)} & : & NE_{1}\otimes NE_{2} &
\longrightarrow & C_{3}NE_{3}/\partial _{4}(NE_{4}\cap D_{4}) \\
&  & \left\{ x_{1}\otimes y_{2}\right\} _{(0)(2,1)} & \longmapsto &
\overline{s_{2}s_{1}x_{1}(s_{1}y_{2}-s_{0}y_{2}-s_{2}y_{2})} \\
\left\{ ~\otimes ~\right\} _{(1)(0)} & : & NE_{2}\otimes NE_{2} &
\longrightarrow & NE_{3}/\partial _{4}(NE_{4}\cap D_{4}) \\
&  & \left\{ x_{2}\otimes y_{2}\right\} _{(1)(0)} & \longmapsto & \overline{%
(s_{1}x_{2}-s_{2}x_{2})s_{1}y_{2}+s_{2}(x_{2}y_{2})} \\
\left\{ ~\otimes ~\right\} _{(2)(0)} & : & NE_{2}\otimes NE_{2} &
\longrightarrow & NE_{3}/\partial _{4}(NE_{4}\cap D_{4}) \\
&  & \left\{ x_{2}\otimes y_{2}\right\} _{(2)(0)} & \longmapsto & \overline{%
-s_{2}x_{2}s_{0}y_{2}} \\
\left\{ ~\otimes ~\right\} _{(2)(1)} & : & NE_{2}\otimes NE_{2} &
\longrightarrow & NE_{3}/\partial _{4}(NE_{4}\cap D_{4}) \\
&  & \left\{ x_{2}\otimes y_{2}\right\} _{(2)(1)} & \longmapsto & \overline{%
s_{2}x_{2}(s_{2}y_{2}-s_{1}y_{2})}%
\end{array}%
\end{equation*}%
(The elements denoted by $\overline{(\text{ },\text{ })}$ are cosets in $%
NE_{3}/\partial _{4}(NE_{4}\cap D_{4})$ and given by the elements in $%
NE_{3}. $)
\end{proposition}

\begin{proof}
Here we will check some of the conditions.The others can be checked easily%
\newline
\newline
\textbf{3CM9)} Since%
\begin{equation*}
\begin{array}{lll}
\partial _{4}\left( C_{\left( 3,2\right) \left( 1,0\right) }(x_{2}\otimes
y_{3}\right) ) & = & \left(
s_{2}s_{0}d_{2}x_{2}-s_{1}s_{0}d_{2}x_{2}-s_{0}x_{2}\right) s_{2}d_{3}y_{3}%
\end{array}%
\end{equation*}
iwe find%
\begin{equation*}
\begin{array}[b]{lll}
\left\{ \overline{\partial }_{2}x_{2}\otimes \overline{\partial }%
_{3}y_{3}\right\} _{(1,0)(2)}^{3} & = & \left(
s_{2}s_{0}d_{2}x_{2}-s_{1}s_{0}d_{2}x_{2}-s_{0}x_{2}\right) s_{2}d_{3}y_{3}
\\
&  & s_{0}x_{2}s_{2}d_{3}y_{3}\in {mod}\partial _{4}(NE_{4}\cap D_{4})
\\
& = & -\left\{ x_{2}\otimes \partial _{3}y_{3}\right\} _{(0)(2)}^{3}%
\end{array}%
\end{equation*}
\newline
\textbf{3CM13) }Since\textbf{\ }%
\begin{equation*}
\begin{array}{lll}
\partial _{4}\left( C_{\left( 3,2\right) \left( 1\right) }(x_{2}\otimes
y_{3}\right)  & = & s_{2}x_{2}(s_{2}y_{3}-s_{1}y_{3}-y_{3})%
\end{array}%
\end{equation*}%
\begin{equation*}
\begin{array}{lll}
& = &
\end{array}%
\end{equation*}%
we find%
\begin{equation}
\begin{array}[b]{lll}
\left\{ x_{2}\otimes \partial _{3}y_{3}\right\} _{(2)(1)}^{3} & = &
s_{2}x_{2}(s_{1}d_{3}y_{3}-s_{2}d_{3}y_{3}) \\
& \equiv  & s_{2}x_{2}y_{3}\in {mod}\partial _{4}(NE_{4}\cap D_{4}) \\
& = & x_{2}\cdot y_{3}%
\end{array}
\tag{3.13}
\end{equation}
\end{proof}

\begin{theorem}
The category of $3$-crossed modules is equivalent to the category of
simplicial algebras with Moore complex of length $3$.
\end{theorem}

\subsection{Projective 3-crossed Resolution}

Here as an application we will define projective $3$-crossed resolution of
commutative algebras. This construction was defined by PJ.L.Doncel, A.R.
Grandjean and M.J.Vale\textsc{\ }in \cite{doncel} for $2$-crossed modules.

\begin{definition}
A projective $\mathit{3}$\textit{-crossed }resolution of an $\mathbf{k}$%
-algebra $E$ is an exact sequence
\begin{equation*}
...\longrightarrow C_{k+1}\overset{\partial _{k}}{\longrightarrow }%
C_{k}\longrightarrow ...\longrightarrow C_{3}\overset{\partial _{3}}{%
\longrightarrow }C_{2}\overset{\partial _{2}}{\longrightarrow }C_{1}\overset{%
\partial _{1}}{\longrightarrow }C_{0}\overset{\partial _{0}}{\longrightarrow
}E\longrightarrow 0
\end{equation*}%
of $k$-modules such that \newline
$1)$ $C_{0}$ is projective in the category of $\mathbf{k}$-algebras\newline
$2)$ $C_{i}$ is a $C_{i-1}$-algebras and projective in the category of $%
C_{i-1}$ algebras for $i=1,2$\newline
$3)$ For any epimorphism $F=(f,id,id,id):(C_{3}^{\prime
},C_{2},C_{1},C_{0},\partial _{3}^{\prime },\partial _{2},\partial
_{1})\longrightarrow (C_{3}^{\prime \prime },C_{2},C_{1},C_{0},\partial
_{3}^{\prime \prime },\partial _{2},\partial _{1})$ and morphism $%
H=(h,id,id,id):(C_{3},C_{2},C_{1},C_{0},\partial _{3},\partial _{2},\partial
_{1})\longrightarrow (C_{3}^{\prime \prime },C_{2},C_{1},C_{0},\partial
_{3}^{\prime },\partial _{2},\partial _{1})$ there exist a morphism $%
Q=(q,id,id,id):(C_{3},C_{2},C_{1},C_{0},\partial _{3},\partial _{2},\partial
_{1})\longrightarrow (C_{3}^{\prime },C_{2},C_{1},C_{0},\partial
_{3}^{\prime },\partial _{2},\partial _{1})$ such that $FQ=H$ \newline
$4)$ for $k\geq 4$, $C_{k}$ is a projective $\mathbf{k}$-module\newline
$5)$ $\partial _{4}$ is a homomorphism of $C_{0}$-module where the action of
$C_{0}$ on $C_{4}$ is defined by $\partial _{0}$\newline
$6)$ For $k\geq 5,$ $\partial _{k}$ is a homomorphism of $k$-modules\newline
\end{definition}

\begin{proposition}
Any commutative $\mathbf{k}$-algebra with a unit has a projective $3$%
-crossed resolution.
\end{proposition}

\begin{proof}
Let $\mathbf{E}$ be a $\mathbf{k}$-algebra and $C_{0}=\mathbf{k}[X_{i}]$ a
polynomial ring such that there exist an epimorphism $\partial
_{0}:C_{0}\rightarrow B$.\newline
Now let define $C_{1}$ as $C_{0}=[ker\partial _{0}]$ the positively graded
part of polynomial ring on $ker\partial _{0}$ and define $\partial
_{1}:C_{1}\rightarrow C_{0}$ by inducing from the inclusion $i:ker\partial
_{0}\rightarrow C_{0}$.\newline
Now let define $K_{2}=C_{0}(C_{1}\times C_{1}\bigcup ker\partial _{1})$, the
free $C_{0}$-module on the disjoint union \newline
$(C_{1}\times C_{1})\bigcup ker\partial _{1}$ and define

\begin{equation*}
\begin{array}{lll}
\partial _{2}^{\prime }:K_{2}\rightarrow ker\partial _{1} &  &
\end{array}%
\end{equation*}%
as%
\begin{equation*}
\begin{array}{lll}
\partial _{2}^{\prime }(x_{1}y_{1})=x_{1}y_{1}-\partial _{2}(y_{1})x_{1},\ \
\ x_{1}y_{1}\in C_{1} &  &  \\
\partial _{2}^{\prime }(x)=x,\ \ \ \ x\in ker\partial _{1} &  &
\end{array}%
\end{equation*}%
Let $\mathbf{R}^{\prime }$ be the $C_{0}$-module generated by the relations
\newline
\begin{equation*}
\begin{array}{lll}
(\alpha x_{1}+\beta y_{1},z_{1})-\alpha (x_{1},z_{1})-\beta (y_{1},z_{1}) &
&  \\
(x_{1},\alpha y_{1}+\beta z_{1})-\alpha (x_{1},y_{1})-\beta (x_{1},z_{1}) &
&
\end{array}%
\end{equation*}%
when $\alpha ,\beta \in C_{0}$ and $x_{1},y_{1},z_{1}\in C_{1}$. Now define $%
C_{2}=K_{2}/R^{\prime }$. Now define $\partial _{2}:C_{2}\rightarrow
ker\partial _{1}$ with $\partial _{2}\pi =\partial _{2}^{\prime }$ where $%
\pi :K_{2}\rightarrow (C_{2}=K_{2}/R^{\prime })$ is projection.\newline
Now we will define $C_{3}$. Let $K_{3}$ is the $C_{0}$-module defined on the
disjoint union \newline
$C_{0}(A_{(1,0)}\cup A_{(0,2)}\cup A_{(2,1)}\cup A_{(1,0)(2)}\cup
A_{(2,0)(1)}\cup A_{(0)(2,1)}\cup ker\partial _{2})$ \newline
where \newline
$A_{(1,0)}=A_{(0,2)}=A_{(2,1)}=C_{2}\times C_{2}$ ,\newline
$A_{(1,0)(2)}=A_{(2,0)(1)}=C_{1}\times C_{2}$ \newline
and $A_{(0)(2,1)}=C_{2}\times C_{1}$.\newline
We have $\partial _{3}^{\prime }:K_{3}\rightarrow (C_{1}\times C_{1}\cup
ker\partial _{1})$
\begin{equation*}
\begin{array}{lll}
\partial _{3}^{\prime }(x,y)=xy-\partial _{2}yx,\ \ \ \ \ (x,y)\in A_{(2)(1)}
&  &  \\
\partial _{3}^{\prime }(x,y)=(\partial _{2}x,\partial _{2}y)+xy,\ \ \ \ \
(x,y)\in A_{(1)(0)} &  &  \\
\partial _{3}^{\prime }(x,y)=(\partial _{2}x,y)+yx,\ \ \ \ \ (x,y)\in
A_{(0)(2,1)} &  &  \\
\partial _{3}^{\prime }(x,y)=(x,\partial _{2}y)-\partial _{1}xy+xy,\ \ \ \ \
(x,y)\in A_{(2,0)(1)} &  &  \\
\partial _{3}^{\prime }(x,y)=0,\ \ \ \ \ (x,y)\in A_{(0)(2)} &  &  \\
\partial _{3}^{\prime }(x,y)=0,\ \ \ \ \ (x,y)\in A_{(1,0)(2)} &  &  \\
\partial _{3}^{\prime }(x,y)=x,\ \ \ \ \ x\in ker\partial _{2} &  &  \\
&  &
\end{array}%
\end{equation*}%
Now define a $C_{0}$-module $R$, generated by the relations,where $\alpha
,\beta \in C_{0}$. Now define $C_{3}=K_{3}/R$ we have $\partial
_{3}:C_{3}\rightarrow ker\partial _{2}$ with $\partial _{3}\pi :\partial
_{3}^{\prime }$ where $\pi :K_{3}\rightarrow (C_{3}=K_{3}/R)$ is projection.
With these constructions
\begin{equation*}
C_{3}\overset{\partial _{3}}{\longrightarrow }C_{2}\overset{\partial _{2}}{%
\longrightarrow }C_{1}\overset{\partial _{1}}{\longrightarrow }C_{0}
\end{equation*}%
is a projective $3$-crossed module. If we define $C_{4}$ as the projection
resolution of the $\mathbf{E}$-module $ker\partial _{3}$ then we have the
projective crossed resolution
\begin{equation*}
\cdots C_{k}{\longrightarrow }C_{k-1}\cdots C_{4}\overset{q}{\longrightarrow
}C_{3}\overset{\partial _{3}}{\longrightarrow }C_{2}\overset{\partial _{2}}{%
\longrightarrow }C_{1}\overset{\partial _{1}}{\longrightarrow }C_{0}\overset{%
\partial _{0}}{\longrightarrow }E{\longrightarrow }0
\end{equation*}%
where $q$ is the projection.
\end{proof}

\subsection{Lie Algebra Case}

Lie algebraic version of crossed modules were introduced by Kassel and Loday
in, \cite{kas}, and that of $2$-crossed modules were introduced by G. J.
Ellis in, \cite{ellis1}. Also higher dimensional Peiffer elements in
simplicial Lie algebras introduced in \cite{aa}. As a consequence of the
commutative algebra version of $3$-crossed modules in this section we will
define the $3$-crossed modules of Lie algebras by using the results given in
\cite{aa} with a similar way used for defining the commutative algebra case.
The relations for commutative algebra case given in the previous section can
be applied to the Lie algebra case with some slight differences in the
proofs up to the definition.

\begin{definition}
A $\mathit{3}$\textit{-crossed module} over Lie algebras consists of a
complex of Lie algebras
\begin{equation*}
L_{3}\overset{\partial _{3}}{\longrightarrow }L_{2}\overset{\partial _{2}}{%
\longrightarrow }L_{1}\overset{\partial _{1}}{\longrightarrow }L_{0}
\end{equation*}%
together with an action of $L_{0}$ on $L_{3},L_{2},L_{1}$ and an action of $%
L_{1}$ on $L_{3},L_{2}$ and an action of $L_{2}$ on $L_{3}$ so that $%
\partial _{3}$, $\partial _{2}$,$\partial _{1}$ are morphisms of $%
L_{0},L_{1} $-groups and the $L_{1},L_{0}$-equivariant liftings
\begin{equation*}
\begin{tabular}{lll}
$\{$ $,$ $\}_{(1)(0)}:L_{2}\times L_{2}\longrightarrow L_{3},$ & $\{$ $,$ $%
\}_{(0)(2)}:L_{2}\times L_{2}\longrightarrow L_{3},$ & $\{$ $,$ $%
\}_{(2)(1)}:L_{2}\times L_{2}\longrightarrow L_{3},$ \\
&  &  \\
$\{$ $,$ $\}_{(1,0)(2)}:L_{1}\times L_{2}\longrightarrow L_{3},$ & $\{$ $,$ $%
\}_{(2,0)(1)}:L_{1}\times L_{2}\longrightarrow L_{3},$ &  \\
&  &  \\
$\{$ $,$ $\}_{(0)(2,1)}:L_{2}\times L_{1}\longrightarrow L_{3},$ & $\{$ $,$ $%
\}:L_{1}\times L_{1}\longrightarrow L_{2}$ &
\end{tabular}%
\end{equation*}%
called \textit{$3$-dimensional Peiffer liftings}. This data must satisfy the
following axioms:%
\begin{equation*}
\begin{array}{lrrl}
\mathbf{3CM1)} &  & C_{3}\overset{\partial _{3}}{\longrightarrow }C_{2}%
\overset{\partial _{2}}{\longrightarrow }C_{1} & \text{is a }2\text{-crossed
module with the Peiffer lifting }\{\text{ }\otimes \text{ }\}_{(2),(1)} \\
\mathbf{3CM2)} &  & \partial _{2}\left\{ l_{1}\otimes m_{1}\right\} & =\text{
}^{\text{ }\partial _{1}m_{1}}l_{1}-\left[ l_{1},m_{1}\right] \\
\mathbf{3CM3)} &  & \left\{ l_{2}\otimes \partial _{2}m_{2}\right\}
_{(0)(2,1)} & =\left\{ l_{2}\otimes m_{2}\right\} _{(2)(1)}-\left\{
l_{2}\otimes m_{2}\right\} _{(1)(0)} \\
\mathbf{3CM4)} &  & \partial _{3}\left\{ l_{2}\otimes m_{2}\right\} _{(1)(0)}
& =\left\{ \partial _{2}l_{2}\otimes \partial _{2}m_{2}\right\} +\left[
l_{2},m_{2}\right] \\
\mathbf{3CM5)} &  & \left\{ l_{1}\otimes \partial _{3}l_{3}\right\}
_{(2,0)(1)} & =\left\{ l_{1}\otimes \partial _{3}l_{3}\right\}
_{(0)(2,1)}+\left\{ l_{1}\otimes \partial _{3}l_{3}\right\} _{(1,0)(2)}-%
\text{ }^{\partial _{1}l_{1}}l_{3} \\
\mathbf{3CM6)} &  & \left\{ \partial _{2}l_{2}\otimes m_{2}\right\}
_{(2,0)(1)} & =-\left\{ l_{2}\otimes m_{2}\right\} _{(0)(2)}+\left[
l_{2},m_{2}\right] \cdot \left\{ l_{2}\otimes m_{2}\right\}
_{(2)(1)}+\left\{ l_{2}\otimes m_{2}\right\} _{(1)(0)} \\
\mathbf{3CM7)} &  & \left\{ \partial _{3}l_{3}\otimes \partial
_{3}m_{3}\right\} _{(1)(0)} & =\left[ m_{3},l_{3}\right] \\
\mathbf{3CM8}) &  & \left\{ \partial _{3}l_{3}\otimes \partial
_{2}l_{2}\right\} _{(0)(2,1)} & =\text{ }^{-\partial _{2}l_{2}}l_{3} \\
\mathbf{3CM9)} &  & \left\{ \partial _{2}l_{2}\otimes \partial
_{3}l_{3}\right\} _{(1,0)(2)} & =-\left\{ l_{2}\otimes \partial
_{3}l_{3}\right\} _{(0)(2)} \\
\mathbf{3CM10)} &  & \left\{ \partial _{2}l_{2}\otimes \partial
_{3}l_{3}\right\} _{(2,0)(1)} & =\text{ }^{\partial _{2}l_{2}}l_{3}-\left\{
l_{2}\otimes \partial _{3}l_{3}\right\} _{(0)(2)} \\
\mathbf{3CM11)} &  & \left\{ \partial _{3}l_{3}\otimes l_{1}\right\}
_{(0)(2,1)} & =\text{ }^{-l_{1}}l_{3} \\
\mathbf{3CM12)} &  & \left\{ l_{2}\otimes \partial _{3}l_{3}\right\}
_{(1)(0)} & =-l_{2}\cdot l_{3} \\
\mathbf{3CM13)} &  & \left\{ \partial _{3}l_{3}\otimes l_{2}\right\}
_{(1)(0)} & =l_{2}\cdot l_{3} \\
\mathbf{3CM14)} &  & \left\{ \partial _{3}l_{3}\otimes l_{2}\right\}
_{(2)(0)} & =0 \\
\mathbf{3CM15)} &  & \partial _{3}\left\{ l_{1}\otimes l_{2}\right\}
_{(2,0)(1)} & =\partial _{3}\left\{ l_{1}\otimes l_{2}\right\}
_{(1,0)(2)}+\left\{ l_{1}\otimes \partial _{2}l_{2}\right\} -\text{ }%
^{\partial _{1}l_{1}}l_{2}+\text{ }^{l_{1}}l_{2} \\
\mathbf{3CM16)} &  & \partial _{3}\left\{ l_{1}\otimes l_{2}\right\}
_{(0)(2,1)} & =\left\{ l_{1}\otimes \partial _{2}l_{2}\right\} -\text{ }%
^{l_{1}}l_{2}%
\end{array}%
\end{equation*}
\end{definition}

T.S. Kuzp\i nar\i   \\            
stufan@aksaray.edu.tr      \\
Department of Mathematics\\ 
         Aksaray University\\        
         Aksaray/Turkey

\bigskip

A. Odaba\c{s} \\
aodabas@ogu.edu.tr      \\
Mathematics and Computer Sciences Department \\
         Eski\c{s}ehir Osmangazi University\\
         Eski\c{s}ehir/Turkey
         
\bigskip

E.\"{O}. Uslu \\
euslu@aku.edu.tr      \\
Department of Mathematics\\
         Afyon Kocatepe University\\
         Afyonkarahisar/Turkey

\end{document}